\documentclass[12pt]{article}

\usepackage{amscd,amsmath, amssymb, fancyhdr, mathbbol}

\numberwithin{equation}{section}

\newcommand{\version}{version 4.0,\ \ Jan. 8, 2020}

\def\eqref#1{(\ref{#1})}

\newcommand{\arrow}{{\:\longrightarrow\:}}
\newcommand{\Z}{{\Bbb Z}}
\def\C{{\Bbb C}}

\newcommand{\Q}{{\Bbb Q}}

\def\1{\sqrt{-1}\:}
\newcommand{\restrict}[1]{{\left|_{{\phantom{|}\!\!}_{#1}}\right.}}
\newcommand{\cntrct}                
{\hspace{2pt}\raisebox{1pt}{\text{$\lrcorner$}}\hspace{2pt}}

\newcommand{\calo}{{\cal O}}

\renewcommand{\phi}{\varphi}
\renewcommand{\epsilon}{\varepsilon}
\renewcommand{\geq}{\geqslant}
\renewcommand{\leq}{\leqslant}

\newcommand{\Pic}{\operatorname{Pic}}
\newcommand{\Ext}{\operatorname{Ext}}
\newcommand{\Hom}{\operatorname{Hom}}

\newcounter{Mycounter}[section]
\newcounter{lemma}[section]
\renewcommand{\thelemma}{{Lemma \thesection.\arabic{lemma}}}
\newcommand{\lemma}{%
    \setcounter{lemma}{\value{Mycounter}}
    \refstepcounter{lemma}
    \stepcounter{Mycounter}
    {\noindent \bf \thelemma:\ }}

\newcounter{claim}[section]
\renewcommand{\theclaim}{{Claim \thesection.\arabic{claim}}}
\newcommand{\claim}{%
    \setcounter{claim}{\value{Mycounter}}
    \refstepcounter{claim}
    \stepcounter{Mycounter}
    {\noindent \bf \theclaim:\ }}

\newcounter{sublemma}[section]
\newcounter{corollary}[section]
\renewcommand{\thecorollary}{{Corollary \thesection.\arabic{corollary}}}
\newcommand{\corollary}{%
    \setcounter{corollary}{\value{Mycounter}}
    \refstepcounter{corollary}
    \stepcounter{Mycounter}
    {\noindent \bf \thecorollary:\ }}

\newcounter{theorem}[section]
\renewcommand{\thetheorem}{{Theorem \thesection.\arabic{theorem}}}
\newcommand{\theorem}{%
    \setcounter{theorem}{\value{Mycounter}}
    \refstepcounter{theorem}
    \stepcounter{Mycounter}
    {\noindent \bf \thetheorem:\ }}

\newcounter{conjecture}[section]
\newcounter{proposition}[section]
\renewcommand{\theproposition}
      {{Proposition \thesection.\arabic{proposition}}}
\newcommand{\proposition}{%
    \setcounter{proposition}{\value{Mycounter}}
    \refstepcounter{proposition}
    \stepcounter{Mycounter}
    {\noindent \bf \theproposition:\ }}

\newcounter{definition}[section]
\renewcommand{\thedefinition}
      {{Definition~\thesection.\arabic{definition}}}
\newcommand{\definition}{%
    \setcounter{definition}{\value{Mycounter}}
    \refstepcounter{definition}
    \stepcounter{Mycounter}
    {\noindent \bf \thedefinition:\ }}

\newcounter{example}[section]
\newcounter{remark}[section]
\renewcommand{\theremark}{{Remark \thesection.\arabic{remark}}}
\newcommand{\remark}{%
    \setcounter{remark}{\value{Mycounter}}
    \refstepcounter{remark}
    \stepcounter{Mycounter}
    {\noindent \bf \theremark:\ }}

\newcounter{problem}[section]
\newcounter{question}[section]
\newcommand{\proof}{\noindent{\bf Proof:\ }}

\makeatletter

\@addtoreset{equation}{section} \@addtoreset{footnote}{section}
\makeatother

\def\blacksquare{\hbox{\vrule width 5pt height 5pt depth 0pt}}
\def\endproof{\blacksquare}

\begin{document}
\begin{center}
{\LARGE\bf
Pullbacks of hyperplane sections for Lagrangian fibrations are primitive
\\[4mm]
}

Ljudmila Kamenova\footnote{Partially 
supported by a grant from the Simons Foundation/SFARI (522730, LK).}, 
Misha Verbitsky\footnote{Partially 
supported by the Russian Academic Excellence Project '5-100'.}

\end{center}
\begin{center}{\em \small
Dedicated to Professor Claire Voisin}
\end{center}
{\small \hspace{0.10\linewidth}
\begin{minipage}[t]{0.85\linewidth}
{\bf Abstract.} 
Let $p:\; M\rightarrow B$ be a Lagrangian fibration on 
a hyperk\"ahler manifold of maximal holonomy (also known as IHS),
and $H$ be the generator of the Picard group of $B$. Assume that $p$
has no multiple fibers in codimension 1. We prove that $p^*(H)$
is a primitive class on $M$. 
\end{minipage}
}

\section{Introduction} 

In this paper we consider a compact
hyperk\"ahler manifold of maximal holonomy admitting a 
holomorphic fibration $\pi:\; M \rightarrow B$. 
The fibration structure is quite restricted due to the work of 
Matsushita, \cite{_Matsushita:fibred_}, who first noticed that the general 
fiber is a Lagrangian abelian variety of half of the dimension of the total 
space. The base has the same rational cohomology as $\C P^n$ and the Picard 
group $\Pic(B)$ has rank one. Assume that multiple fibers of $\pi$ have codimension $\geq 2$ in $B$. We prove that the pullback of 
the fundamental class of a hyperplane section is primitive, i.e., 
indivisible as an integral class. 

\hfill

\theorem 
Let $M$ be a hyperk\"ahler manifold admitting a Lagrangian fibration 
$\pi:\; M \rightarrow B$ and $H$ be the generator of $\Pic(B)$.  
Assume that the set of points $b \in B$ with multiple fibers $\pi^{-1}(b)$ 
has codimension at least $2$. 
Then the class $\pi^*H \in H^2(M, \Z)$ is primitive. 

\hfill

The proof is based on the observation that if $\pi^*H \in H^2(M, \Z)$ is 
not primitive, i.e., $\pi^*H = mH'$, then $H'$ has trivial 
cohomology by Demailly, Peternell and Schneider's theorem. 
The assumption excluding multiple fibers is needed to exclude the case 
when $H'$ is trivial on the generic fiber. 
Applying the Hirzebruch-Riemann-Roch formula for an irreducible 
hyperk\"ahler manifold, one would obtain a contradiction.

\section{Basic notions}

\definition
A {\em hyperk\"ahler manifold of maximal holonomy} 
(or {\em irreducible holomorphic symplectic}) 
manifold 
$M$ is a compact complex simply connected K\"ahler manifold with 
$H^{2,0}(M)=\C \sigma$, where $\sigma$ is everywhere non-degenerate. 

\hfill

For the rest of the paper we consider hyperk\"ahler manifolds 
of maximal holonomy. 
Due to the work of Matsushita we know that the fibration structure of 
hyperk\"ahler manifolds is quite restricted. 

\hfill

\begin{theorem} (D. Matsushita, \cite{_Matsushita:fibred_}) \label{fibr}
Let $M$ be a hyperk\"ahler manifold and $f\colon M\rightarrow B$ a 
proper surjective morphism with a smooth base $B$. Assume that $f$ has 
connected fibers and $0 < \dim B < \dim M$. Then $f$ is Lagrangian 
and $\dim_{\C} B = n$, where $\dim_{\C} M = 2n$. 
\end{theorem}

\hfill

\definition
Following \ref{fibr}, we call the morphism $f\colon M\rightarrow B$ a 
{\em Lagrangian fibration} on the hyperk\"ahler manifold $M$. 

\hfill

\remark\label{_CP^n_fibra_Remark_}
In \cite{_Matsushita:CP^n_}, D. Matsushita also proved that 
the base $B$ of a Lagrangian fibration has the same (rational)
cohomology as $\C P^n$. In \cite{_Hwang:CP^n_},
J.-M. Hwang proved that when $B$ is smooth, then it is actually isomorphic
to $\C P^n$.

\hfill

\definition
Given a hyperk\"ahler manifold $M$, there is a non-degenerate primitive 
form $q$ on $H^2(M,\Z)$, called the {\em Beauville-Bogomolov-Fujiki form} 
(or {\em ``BBF form''} for short) 
of signature $(3,b_2-3)$, satisfying the {\em Fujiki relation} 
$$\int_M \alpha^{2n}=c\cdot q(\alpha)^n\qquad\text{for }\alpha \in 
H^2(M,\Z),$$ with $c>0$ a constant depending on the topological type of $M$. 
This form generalizes the intersection pairing on K3 surfaces. 
A detailed description of the form can be found in 
\cite{_Beauville_}, \cite{_Bogomolov:defo_} and \cite{_Fujiki:HK_}. 

\hfill

\definition
Let $f:\; M \arrow B$ be a Lagrangian fibration. 
As shown in \cite{_Matsushita:CP^n_}, $H^*(B,\Q)= H^*(\C P^n, \Q)$.
Let $H$ be a primitive integer generator of $H^2(M, \Q)$,
and $\calo(1)$ be a holomorphic line bundle on $B$
with first Chern class $H$. When $B=\C P^n$ (this is the case when $B$ 
is smooth by Hwang's result \cite{_Hwang:CP^n_}),
the bundle $\calo(1)$ coincides with the usual $\calo(1)$.
We call $H$ {\em the fundamental class of a hyperplane section}.

\hfill

\remark
A {\em semiample bundle} is a base point free line bundle
which has positive Kodaira dimension.
Let $f:\; M \rightarrow B$ be a Lagrangian fibration, and
$L:=f^*(\calo(k))$, $k>0$. Clearly, $L$ is a semiample nef line bundle.
By Matsushita's theorem any semiample nef line bundle is either ample
or obtained this way. The {\em SYZ conjecture}
(due to Tyurin, Bogomolov, Hassett, Tschinkel, Huybrechts and Sawon; see
\cite{_Verbitsky:SYZ_}) claims that the converse is also true:
any nef line bundle on a hyperk\"ahler manifold is either ample 
or semiample. This conjecture is a special case
of Kawamata's abundance conjecture.

\hfill

\remark \label{HRR}
The Hirzebruch-Riemann-Roch formula for an irreducible 
hyperk\"ahler manifold $M$ states that for a line bundle $L$ on $M$, 
$\chi (L) = \sum a_i q(c_1(L))^i$, where the coefficients $a_i$ are 
constants depending on the topology of $M$ (see \cite[Section 1.11]{Huy}). 
In particular, if $q(c_1(L))=0$, then $\chi(L) = a_0 = \chi({\cal O}_M) 
= \sum (-1)^i h^{0,i}(M) = n+1$, where $2n = \dim M$ 
(see \cite[Section 1.7]{Huy}). 

\hfill

Using the Hirzebruch-Riemann-Roch formula, we can easily obtain our main 
result for K3 surfaces. We are grateful to Claire Voisin for this 
observation.

\hfill

\lemma 
Let $S$ be a K3 surface with an elliptic fibration 
$\pi:\; S \rightarrow \C P^1$. Then the class 
$\pi^* \calo(1) \in H^2(S, \Z)$ is primitive. 

\hfill

Indeed, if we assume that $\pi^* \calo(1) = m H$ for $m > 1$, then 
$H$ would be $m$-torsion on all fibers and $h^0(H)=0$. By Serre 
duality, $h^2(H)=0$. Then applying the Hirzebruch-Riemann-Roch formula 
as in \ref{HRR}, we obtain $2 = \chi (H) = h^0(H) - h^1(H) + h^2(H) = 
-h^1(H) \leq 0$ - a contradiction. 

\hfill

\remark
This result was proven by E. Markman for manifolds of K3${}^{[n]}$-type
(\cite{_Markman:fibrations_K3^n_}) and by B. Wieneck for generalized Kummer varieties
(Lemma 2.7 in \cite{W}). We thank Klaus Hulek for pointing out 
this reference to us.

\hfill

Here we restate a theorem by Demailly, Peternell and Schneider applied to 
compact K\"ahler manifolds with trivial canonical bundle, which is the 
set-up we need. The more general version of this result 
is Theorem 2.1.1 in \cite{_DPS_}.
This theorem was obtained under various hypotheses during the 1990s, see
\cite{_Enoki:semipositive_} and
\cite{_Mourougane_}. It was proved in \cite{_Takegoshi_} when $L$ is nef.

\hfill

\theorem (\cite[Corollary 2.1.2]{_DPS_}) \label{DPS} 
Let $(M, I, \omega)$ be a compact K\"ahler manifold,
$K_M$ its canonical bundle, $\dim_{\C} M = N$,
and $E$ a non-trivial nef line bundle on $M$. Assume that $E$ admits a 
Hermitian metric with semipositive curvature form. 
The cohomology class $\omega$ is considered as an element in $H^1(\Omega^1 M).$
Consider the corresponding multiplication operator $\eta \arrow \omega \wedge \eta$ 
maping $H^p(\Omega^{q} M\otimes E)$ to $H^{p+1}(\Omega^{q+1} M\otimes E)$. Then 
$\eta \arrow \omega^i \wedge \eta$ induces a 
surjective map
$H^0(\Omega^{N-i} M \otimes E) \twoheadrightarrow H^i (E\otimes K_M)$.

\hfill

Further on, we shall also need the following trivial topological observation.

\hfill

\claim\label{_H^2_torsion_free_Claim_}
Let $M$ be a hyperk\"ahler manifold of maximal holonomy.
Then $H^2(M)$ is torsion-free.

\hfill

{\bf Proof:} The universal coefficients formula gives
the exact sequence: 
\[
  0 \to \Ext_\Z^1(H_1(X; \Z), \Z) \to H^2(X; \Z) \to
    \Hom_\Z(H_2(X; \Z), \Z)\to 0.
\]
Since $H_1(X, \Z)=0$ for a maximal holonomy hyperk\"ahler manifold, 
this gives an isomorphism $H^2(X;\Z)=\Hom_\Z(H_2(X; \Z), \Z)$,
hence the torsion vanishes.
\endproof


\section{Main Results}


An ``abelian fibration'' below is a holomorphic map with general
fiber an abelian variety (or, in fact, a compact torus).

\hfill

\proposition\label{_restriction_non-zero_Proposition_}
Let $M$ be a smooth manifold admitting an abelian fibration 
$\pi:\; M \rightarrow B$, and $E$ a line bundle on $M$ which is
trivial on generic fibers
and torsion on all fibers of $\pi$. 
Assume that multiple fibers of $\pi$ have codimension $\geq 2$ in $B$. 
Then $E=\pi^* E'$, where $E'$ is a line bundle
on the base.

\hfill

\proof
We need to show that each point $b\in B$ has a neighbourhood $U_b$ such that
the restriction of $E$ to $\pi^{-1}(U_b)$ is a trivial bundle. 

In the sequel, we 
think of the restriction of $E$ to $\pi^{-1}(U_b)$ as of a holomorphic bundle underlying a local
system. This is done by introducing an appropriate flat metric and
proving that the monodromy of its Chern connection is trivial.

The bundle $E$ is torsion on all fibers of $\pi$, hence its tensor
power $E^{\otimes k}$ is trivial on all fibers. Therefore,
 $E^{\otimes k}= \pi^*\pi_*(E^{\otimes k})$ belongs to $\pi^*(\Pic(B))$. 
Given $b\in B$, choose a neighbourhood $U_b\subset B$
such that the restriction of $E^{\otimes k}$ to $\pi^{-1}(U_b)$ is trivial.
To finish the proof it suffices to show that the restriction of $E$ to $W:=\pi^{-1}(U_b)$ is trivial.

Choose a constant metric $h^k$ on $E^{\otimes k}\restrict W=\calo_W$
and let $h$ be its $k$-th root, which is a metric on $E\restrict W$. Since  $h^k$ is 
constant, its curvature is flat, and
the Chern connection $\nabla$ associated with $h$ is also flat.

To finish the proof, it remains to show that monodromy of $\nabla$ is trivial
on all fibers. However, since $E$ is trivial on generic fibers, the
monodromy of $E$ is trivial on generic fibers.

Any special fiber $F_s$ of $\pi$ is a deformation
retract of its neighbourhood $U_s$ (see \cite{_Morrison:Clemens_Schmid_},
\cite{_Persson:degene_}, \cite{_Clemens:degene_}). 
This retraction gives a map from the fundamental
group of a general fiber to the fundamental group of the special fiber:
the general fiber, denoted $F_g$, is embedded to $U_s$, which is
then retracted to $F_s$. This map is clearly surjective if $F_s$ is a simple fiber.
Therefore, the monodromy representation of $E$ on $F_s$ is induced by the 
monodromy of $E$ on $F_g$, which is trivial. We have shown that $E\restrict W=\calo_W$
in a neighbourhood of a simple fiber. 

To deal with multiple fibers, we notice that they occur only in 
a subset $B_0\subset B$ of codimension $>1$.
Outside of the set $B_0$ of multiple fibers, $E$ is the pullback of a bundle
$E'$ on $B \backslash B_0$. However, for any line bundle $L$, and for 
any subvariety $Z\subset M$ of codimension $\geq 2$, 
one has $L= j_*j^* L$, where $j:\; M \backslash Z \arrow M$
is the open embedding (this is also called ``Serre's condition S2'', 
see \cite[Ch. II, Lemma 1.1.12]{_OSS_}). With this observation, we see 
that $E$ is the pullback of a line bundle on the base $B$. 
\endproof

\hfill

\theorem \label{_primitive_theorem_}
Let $M$ be a hyperk\"ahler manifold admitting a Lagrangian fibration 
$\pi:\; M \rightarrow B$, and $H$ the generator of $\Pic(B)$ 
(this group has rank 1, as shown by D. Matsushita; 
see \ref{_CP^n_fibra_Remark_}). 
Assume that multiple fibers of $\pi$ have codimension $\geq 2$ in $B$. 
Then the class $\pi^*H \in H^2(M, \Z)$ is primitive.

\hfill

\begin{proof}
Suppose that $\pi^*H$ is not primitive, and $\pi^*H=kH'$ in $H^2(M, \Z)$.
Denote by $E$ the line bundle with $c_1(E)=H'$. 
By \ref{_restriction_non-zero_Proposition_}, 
$E= \pi^* E'$ unless $E$ is a non-trivial
torsion bundle on general fibers of $\pi$.
In the first case, $\pi^*H$ is primitive. 
This implies that
$E$ is a non-trivial torsion bundle on smooth
fibers of $\pi$, and torsion on all irreducible
components of non-smooth fibers. In this case $c_1(E')=H$, which 
implies that $\pi^*H$ is primitive. 

Let us apply the
Enoki-Mourugane-Takegoshi-Demailly-Peternell-\-Schnei\-der
vanishing theorem
(\ref{DPS}) to the manifold $M$ with the torsion nef line 
bundle $E$ to obtain the surjective map 
$H^0(\Omega^{2n-i} M \otimes E) \twoheadrightarrow H^i (E)$. 
The bundle $TM$ restricted to a regular fiber $S$ of $\pi:\;M \arrow B$  
can be expressed as an extension
\[
0\arrow TS \arrow TM\restrict S \arrow N_S\arrow 0,
\]
where $N_S$ is the normal bundle, which is trivial 
because $S$ is a fiber of the submersion $\pi:\;M \arrow B$, 
in $s\in B$ which gives $N_S=\pi^*T_s B$. However, 
$TS$ is dual to $N_S$, because $S$ is Lagrangian, 
hence $\Omega^{1} M\restrict S$ is an extension of trivial bundles.
A tensor power of a trivial bundle is trivial, and therefore
$\Omega^kM \restrict S$ is also an extension of trivial bundles.
Then $\Omega^{k}(M)\restrict S\otimes E$ has no sections for all $k$, and 
\ref{DPS} implies that $H^i(E)=0$ for all $i$.

To finish the proof, we apply the Hirzebruch-Riemann-Roch formula for the 
hyperk\"ahler manifold $M$ with the line bundle $E$. Since 
$q(E) = q(L) = 0$ and $H^i(E)=0$, from \ref{HRR} we obtain 
$n+1 = \chi(E) = \sum (-1)^i \dim H^i(E) = 0$, a contradiction. 
Therefore, $\pi^* H$ is primitive. 
\end{proof}

\section{Applications}

In this section we describe some applications of the primitivity result. 

\hfill

\proposition \label{pi2_lemma}
Let $M$ be a hyperk\"ahler manifold admitting a Lagrangian fibration 
$f\colon M \rightarrow \C P^n$. Assume that multiple fibers of $f$ have 
codimension $\geq 2$ in $\C P^n$. Then the map $\pi_2(M)  \arrow \pi_2(\C P^n)$ 
is surjective. 

\hfill

{\bf Proof:}
Since we are in the settings of \ref{_primitive_theorem_}, we know that 
$L = f^* {\cal O}(1)$ is primitive, 
i.e., $c_1(L)$ is not divisible. By Poincar\'e duality there is 
$\alpha \in H_2(M, \Z)$ such that the pairing 
$\langle c_1(L), \alpha \rangle = 1$ in $M$. 
This is the same as the pairing 
$\langle f_* \alpha, c_1({\cal O}(1)) \rangle$ in $\C P^n$, 
which means that $f_* \alpha$ is the class of a line, therefore 
$H_2(M, \Z) \arrow H_2(\C P^n, \Z)$ is surjective. 
Since $M$ and $\C P^n$ are simply connected, this induces a surjection on the 
homotopy groups $\pi_2(M) \arrow \pi_2(\C P^n)$ 
(see \cite[Corollary 10.8]{_Bredon_}). 
\endproof

\hfill

\remark
We conjecture that if the fibration $f\colon M \rightarrow \C P^n$ has no 
multiple fibers, then for a general curve $C \subset \C P^n$ there is a 
continuous section $C \rightarrow M$. The evidence is that for every curve 
class $[C]$ there is a class on $M$ surjecting to $[C]$ by \ref{pi2_lemma}. 

\hfill

\definition
A pullback of a very ample bundle is called
{\em very semiample}.

\hfill

\corollary
Let $E$ be a semiample line bundle on a hyperk\"ahler manifold,
which is not ample. Assume that the corresponding
Lagrangian fibration has base $\C P^n$ and that the set of multiple fibers 
has codimension $\geq 2$ in $\C P^n$. Then $E$ is very semiample.

\hfill

\proof 
Indeed, by \ref{_primitive_theorem_}, $E=f^* \calo(i)$, where $i>0$ and
$f:\; M \arrow \C P^n$ is a Lagrangian fibration.
\endproof

\hfill

\claim \label{intro_app}
Let $\pi_1$, $\pi_2$ be 
different Lagrangian fibrations on $M$ with base $\C P^n$, such that the sets of 
multiple fibers of both fibrations have codimension $\geq 2$ in $\C P^n$,
and $P_1, P_2$ be non-trivial nef bundles which are trivial on the 
fibers of $\pi_1, \pi_2$. Then the natural map $F:\; M\arrow {\Bbb P}\big(H^0(M,P_1 \otimes P_2)^*\big)$ is
holomorphic and birational to its image. 

\hfill

\proof 
By \ref{_primitive_theorem_},  $P_i=\pi_i^* \calo(k_i)$, for $i=1,2$. 
Since $P_i$ are nef and non-trivial, one has $k_i>0$, and the line bundle
$P_1 \otimes P_2$ is base point free.
Therefore, the corresponding map $F:\; M\arrow {\Bbb P}\big(H^0(M,P_1 \otimes P_2)^*\big)$
is holomorphic. This map contracts only subvarieties $Z\subset M$ which lie
in fibers of $\pi_1$ and $\pi_2$. Since the corresponding Lagrangian fibers
are transversal when smooth, a general point is not contained in such a subvariety.
\endproof

\hfill

{\bf Acknowledgments.} We are very grateful to Claire Voisin whose 
ideas inspired the proof of the main theorem. The work was 
completed at the SCGP during the second-named author's visit. 
We are grateful to the SCGP for the hospitality. 
The first named author thanks Michel Brion for their conversations 
about an earlier argument of the main theorem and for his interest. 
We would like to thank Christian Lehn for his interest and references.
We are grateful to Ulrike Rie\ss\  for invaluable advice, her interest
to this work, and for finding inconsistencies in an earlier
version of our arguments. We thank the referee for their detailed comments and 
corrections. Much gratitude to D. Kaledin and D. Huybrechts for
finding errors in the arxiv version and their suggestions.

{\small

\noindent {\sc Ljudmila Kamenova\\
Department of Mathematics, 3-115 \\
Stony Brook University \\
Stony Brook, NY 11794-3651, USA,} \\
\tt kamenova@math.sunysb.edu
\\

\noindent {\sc Misha Verbitsky\\
{\sc Universit\'e libre de Bruxelles, CP 213,\\
Bd du Triomphe, 1050 Brussels, Belgium}, \\ {\em also:}\\
{\sc Laboratory of Algebraic Geometry,\\
National Research University HSE,\\
Faculty of Mathematics, 7 Vavilova Str., \\Moscow, Russian Federation,}\\
\tt verbit@verbit.ru
}}

\end{document}